\documentclass[11pt]{article}

\usepackage{tikz}
\usepackage{subfigure}
\usepackage{scalefnt}
\usepackage{amsmath,amsfonts,amsthm}
\usepackage{enumerate}
\usepackage{hyperref}

\numberwithin{equation}{section}

\usepackage{authblk}

\newcommand{\FF}{ \ensuremath{\mathbb{F}}}

\newcommand{\ZZ}{ \ensuremath{\mathbb{Z}}}

\newcommand{\TPSS}{S^{\hspace{.2mm}3} \mbox{$\times
\hspace{-2.5mm}_{-}$} \, S^{\hspace{.1mm}1}}

\newcommand{\TPSSD}{S^{\hspace{.2mm}d-1} \mbox{$\times
\hspace{-2.5mm}_{-}$} \, S^{\hspace{.1mm}1}}

\newcommand\simtimes{\mathbin{%
    \stackrel{\smash{\times}\rule{0pt}{0.9ex}}{\_}%
 }}

\begin{document}

\title{Tight and stacked triangulations of manifolds}
\author{Basudeb Datta}
\affil{Department of Mathematics, Indian Institute of Science, Bangalore 560\,012,
India. dattab@math.iisc.ernet.in.}

\date{}

\maketitle

\vspace{-15mm}

\begin{center}

\date{June 01, 2015}

\end{center}


\begin{abstract}

Tight triangulated manifolds are generalisations of neighborly triangulations of closed surfaces  and are interesting objects in Combinatorial Topology. Tight triangulated manifolds are conjectured to be minimal. Except few, all the known tight triangulated manifolds are  stacked. It is known that locally stacked tight triangulated manifolds are strongly minimal. Except for three infinite series and neighborly surfaces, very few tight triangulated manifolds are known. From some recent works, we know more on tight triangulation. In this article, we present a survey on the works done on tight triangulation. In Section 2, we state some known results on tight triangulations. In Section 3, we present all the known tight triangulated manifolds. Details are available in the references mentioned there. In Section 1, we present some essential definitions. 
\end{abstract}

\section{Preliminaries} 

If $X$ is a $d$-dimensional simplicial complex then, for $0\leq j \leq d$, the number of its $j$-faces is denoted by $f_j = f_j(X)$. The number 
$\chi(X) := \sum_{i=0}^{d} (-1)^i f_i$ is called the \textit{Euler characteristic} of $X$. A simplicial complex $X$ is said to be 
\textit{$k$-neighborly} if $f_{k-1}(X) = \binom{f_0(X)}{k}$. If $X$ is 
2-neighborly then we say $X$ is \textit{neighborly}. 

For any field $\FF$, a closed triangulated $d$-manifold $X$ is said to be
\textit{$\FF$-orientable} if $H_d(X; \FF) \neq \{0\}$.

A simplicial complex $X$ is called \textit{minimal} if $f_0(X) \leq f_0(Y)$, for all triangulation $Y$ of the topological space $|X|$. A simplicial complex $X$ is called \textit{strongly minimal} if $f_i(X) \leq f_i(Y)$, for $0\leq i \leq \dim(X)$ and for all triangulation $Y$ of the topological space $|X|$.

A triangulated $(d+1)$-manifold with non-empty boundary is said to be \textit{stacked} if all its interior faces
have dimension $\geq d$. A closed triangulated $d$-manifold $M$ is said to be \textit{stacked} if $M = \partial
N$ for some stacked triangulated $(d+1)$-manifold $N$. By a \textit{stacked sphere} we mean a stacked
triangulated sphere. A triangulated manifold is said to be \textit{locally stacked} if each vertex link is a
stacked  sphere. The class $\mathcal{K}(d)$ is the set of all locally stacked closed triangulated $d$-manifolds.
Let $\overline{{\mathcal K}}(d)$ be the class of all simplicial complexes whose vertex-links are stacked
$(d-1)$-balls. So, a member of $\overline{{\mathcal K}}(d)$ is a triangulated $d$-manifold with boundary.

Let $X$ be a closed connected triangulated $d$-manifold. Let $\sigma_1$, $\sigma_2$ be $d$-faces and $\psi :
\sigma_1 \to \sigma_2$ a bijection such that $u$ and $\psi(u)$ have no common neighbor in $X$ for each $u\in
\sigma_1$. The simplicial complex $X^{\psi}$ obtained from $X\setminus\{\sigma_1,\sigma_2\}$ by identifying $u$
with $\psi(u)$ for all $u\in \sigma_1$ is a triangulated $d$-manifold and is called the manifold obtained from
$X$ by a \textit{combinatorial handle addition}. For $d\geq 3$, we recursively define the class
$\mathcal{H}^d(k)$ as follows. (a) $\mathcal H^d(0)$ is the set of stacked $(d-1)$-spheres. (b) A triangulated
$d$-manifold $Y$ is in $\mathcal H^d(k +1)$ if it is obtained from a member of $\mathcal H^d(k)$ by a
combinatorial handle addition. (c) The \textit{Walkup's class} {$\mathcal H^d$} is the union $\mathcal
H^d=\bigcup_{k \geq 0} \mathcal H^d(k)$.

For a field $\FF$, a simplicial complex $X$ is called \textit{$\FF$-tight} if (i) $X$ is connected, and (ii) for
all induced subcomplexes $Y$ of $X$ and for all $0\leq j \leq \dim(X)$, the morphism $H_{j}(Y; \mathbb{F}) \to
H_{j}(X; \FF)$ induced by the inclusion map $Y \hookrightarrow X$ is injective.

For a triangulated $d$-manifold $M$, let $g_2(M) := f_1(M) - (d+1)f_0(M) + \binom{d+2}{2}$. So, $g_2(M) \leq
\binom{f_0(M)}{2} - (d+1)f_0(M) + \binom{d+2}{2} = \binom{~f_0(M)-d- 1~}{2}$. If $M$ is $\FF$-orientable manifold
of dimension $d\geq 3$ then $g_2(M) \geq \binom{d+2}{2}\beta_1(M; \FF)$ \cite{ns09}. For $d\geq 3$, a
triangulated $d$-manifold $M$ is called \textit{tight-neighborly} if $\binom{~f_0(M)-d- 1~}{2} =\binom{d
+2}{2}\beta_1(M; \mathbb{Z}/2\ZZ)$. So, $M$ is tight-neighborly if and only if $M$ is neighborly and $g_2(M) =
\binom{d+2}{2}\beta_1(M; \FF)$.

\section{Some known results}

\begin{description}
\item[R01.] Let $M$ be a closed triangulated manifold of dimension $d\geq 2$. If $M$ is $\FF$-tight then $M$ is
$\FF$-orientable and neighborly \cite{bd16}.

\item[R02.] A closed triangulated 2-manifold $M$ is $\FF$-tight $\Longleftrightarrow$ $M$ is $\FF$-orientable and
neighborly \cite{bd16}.

\item[R03.] Let $M$ be a closed triangulated 3-manifold. Then
\begin{align*}
M\in \mathcal{H}^4 \stackrel{\cite{dm14}}{\Longleftrightarrow} M \mbox{ is stacked }
\stackrel{\cite{ba14}}{\Longleftrightarrow} g_2(M) = 10\beta_1(M; 
\ZZ/2\ZZ).
\end{align*}

\item[R04.] Let $M$ be a locally stacked, $\FF$-orientable, closed manifold of dimension $d\geq 3$. Then $M$ is
tight-neighborly $\Longleftrightarrow$ $M$ is $\FF$-tight  \cite{bd16}.

\item[R05.] Let $M$ be a closed  triangulated 3-manifold. Then
\begin{enumerate}
\item[\mbox{(a)}] $M$ is neighborly member of $\mathcal{K}(3)$ and is $\FF$-orientable $\not\Longrightarrow M$ is
$\FF$-tight \cite{bd15, bds1}.

\item[\mbox{(b)}] $M$ is neighborly member of $\mathcal{K}(3) \not\Longrightarrow M$ is stacked \cite{bds1}.
\end{enumerate}

\item[R06.] Let $M$ be a closed triangulated 3-manifold and $\FF$ a field. Then
\begin{enumerate}
\item[\mbox{(a)}] $M$ is tight-neighborly $\stackrel{\cite{bdss}}{\Longrightarrow} M \in \mathcal{K}(3)$ and $M$ is $(\ZZ/2\ZZ)$-tight.

\item[\mbox{(b)}] $M$ is $\FF$-tight, ${\rm char}(\FF) \neq 2 \stackrel{\cite{bds1}}{\Longrightarrow} M\in
\mathcal{K}(3)$.
\end{enumerate}


\item[R07.] Let $M$ be a closed, $\FF$-orientable triangulated 3-manifold.
\begin{enumerate}
\item[\mbox{(a)}] $M$ is neighborly member of $\mathcal{H}^4$ $\stackrel{\cite{dm14}}{\Longleftrightarrow} M$ is stacked and neighborly
$\stackrel{\cite{ba14}}{\Longleftrightarrow} M$ is tight-neighborly $\stackrel{\cite{bdss, bd16}}{\Longleftarrow\!\!\Longrightarrow} M$ is $\FF$-tight and in $\mathcal{K}(3)$.

\item[\mbox{(b)}] If ${\rm char}(\FF) \neq 2$ then 
$M$ is tight-neighborly $\stackrel{\cite{bdss, bd16, bds1}}{\Longleftarrow\!=\!=\!\Longrightarrow} M$ is $\FF$-tight.
\end{enumerate}

\item[R08.] Let $M$ be a closed triangulated manifold of dimension $d\geq 4$. Then


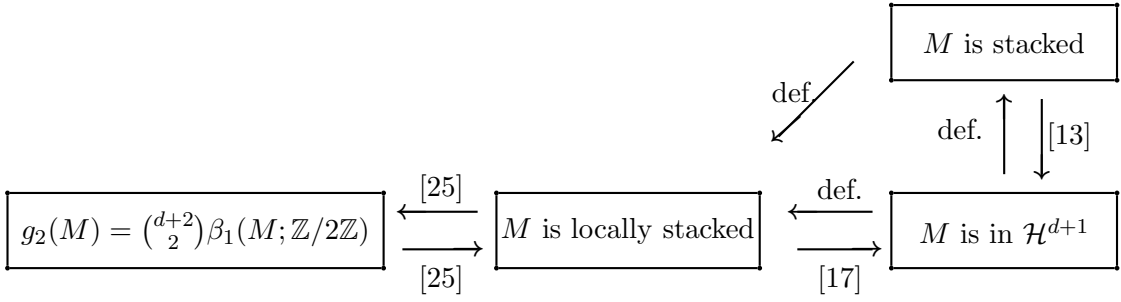
\begin{figure}[htbp]

\tikzstyle{ver}=[] \tikzstyle{extra}=[circle, draw, fill=black!50, inner sep=0pt, minimum width=0pt]
 \tikzstyle{edge} = [draw,thick,-]
 \centering

 \begin{tikzpicture}[black,thick,scale=0.5]

 \begin{scope}[shift={(-23,0)}]
\foreach \x/\y/\z in {5/-1/1,5/1/2,-5/1/3,-5/-1/4}{ \node[extra] (\z) at (\x,\y){};} \foreach \x/\y in
{1/2,2/3,3/4,4/1}{ \path[edge] (\x) -- (\y);}

\node[ver] () at (0,0){\mbox{$g_2(M) = \binom{d+2}{2}\beta_1(M; 
\ZZ/2\ZZ)$}};
\end{scope}

\begin{scope}[shift={(-1.5,5)}]
\foreach \x/\y/\z in {3/-1/1,3/1/2,-3/1/3,-3/-1/4}{ \node[extra] (\z) at (\x,\y){};} \foreach \x/\y in
{1/2,2/3,3/4,4/1}{ \path[edge] (\x) -- (\y);}

\node[ver] () at (0,0){\mbox{{$M$} is stacked}};

\end{scope}

\begin{scope}[shift={(-11.5,0)}]
\foreach \x/\y/\z in {3.5/-1/1,3.5/1/2,-3.5/1/3,-3.5/-1/4}{ \node[extra] (\z) at (\x,\y){};} \foreach \x/\y in
{1/2,2/3,3/4,4/1}{ \path[edge] (\x) -- (\y);}

\node[ver] () at (0,0){\mbox{{$M$} is locally stacked}};

\end{scope}

\begin{scope}[shift={(-1.5,0)}]
\foreach \x/\y/\z in {3/-1/1,3/1/2,-3/1/3,-3/-1/4}{ \node[extra] (\z) at (\x,\y){};} \foreach \x/\y in
{1/2,2/3,3/4,4/1}{ \path[edge] (\x) -- (\y);}

\node[ver] () at (0,0){\mbox{{$M$} is in {$\mathcal{H}^{d+1}$}}};

\end{scope}

\begin{scope}[]




\path[edge] (-7.5,2.5) -- (-5.5,4.5);

\node[ver] () at (-7.2,3.7){\mbox{ def.}};

\node[ver] () at (-7.3,2.7){\mbox{\boldmath{$\swarrow$}}};

\path[edge] (-7,0.5) -- (-5,0.5); \path[edge] (-7,-0.5) -- (-5,-0.5);

\path[edge] (-15.5,0.5) -- (-17.5,0.5); \path[edge] (-17.5,-0.5) -- (-15.5,-0.5);


\path[edge] (-1.5,1.5) -- (-1.5,3.5); \path[edge] (-0.5,1.5) -- (-0.5,3.5);

\node[ver] () at (-6.85,0.47){\mbox{\boldmath{$\leftarrow$}}};

\node[ver] () at (-5.15,-0.53){\mbox{\boldmath{$\rightarrow$}}};
 \node[ver] () at (-17.25,0.47){\mbox{\boldmath{$\leftarrow$}}};
 \node[ver] () at (-15.75,-0.53){\mbox{\boldmath{$\rightarrow$}}};
 \node[ver] () at (-1.5,3.3){\mbox{\boldmath{$\uparrow$}}};
 \node[ver] () at  (-0.5,1.7){\mbox{\boldmath{$\downarrow$}}};
 \node[ver] () at (-6,1.1){\mbox{ def.}};
 \node[ver] () at (-6,-1.3){\mbox{ {$\cite{ka87}$}}};
 \node[ver] () at (-16.6,1.1){\mbox{ {$\cite{ns09}$}}};
 \node[ver] () at (-16.6,-1.3){\mbox{ {$\cite{ns09}$}}};
 \node[ver] () at (-2.8,2.7){\mbox{ def.}};
 \node[ver] () at (0.1,2.5){\mbox{ {$\cite{dm14}$}}};
\end{scope}
\end{tikzpicture}

\end{figure}

\newpage 

\item[R09.] Let $M$ be a closed $\FF$-orientable triangulated manifold of dimension $d\geq 4$. Then


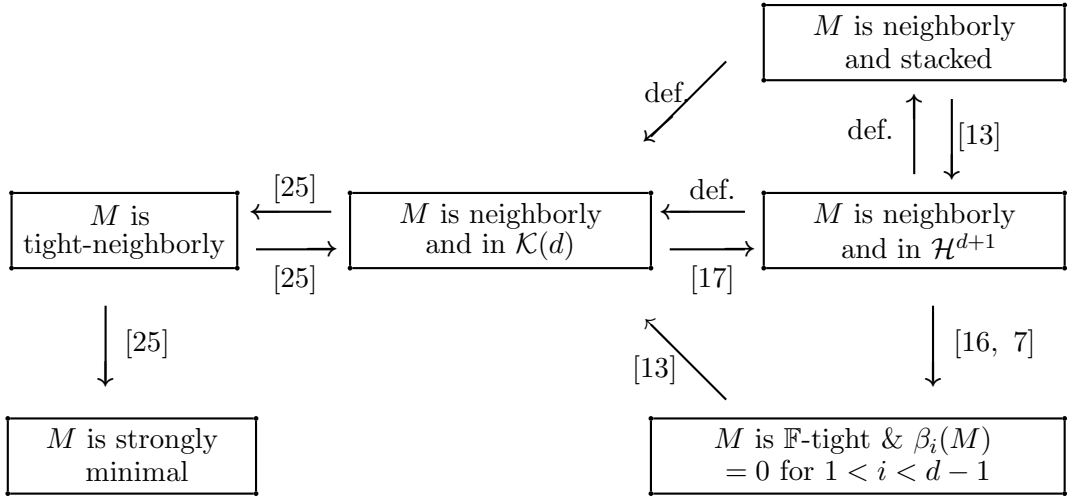
\begin{figure}[htbp]

\tikzstyle{ver}=[] \tikzstyle{extra}=[circle, draw, fill=black!50, inner sep=0pt, minimum width=0pt]
 \tikzstyle{edge} = [draw,thick,-]
 \centering

 \begin{tikzpicture}[black,thick,scale=0.5]

\begin{scope}[shift={(-0.5,5)}]
\foreach \x/\y/\z in {4/-1/1,4/1/2,-4/1/3,-4/-1/4}{ \node[extra] (\z) at (\x,\y){};} \foreach \x/\y in
{1/2,2/3,3/4,4/1}{ \path[edge] (\x) -- (\y);}

\node[ver] () at (0,0.4){\mbox{{$M$} is neighborly}};

\node[ver] () at (0,-0.4){\mbox{ and stacked}};
\end{scope}

\begin{scope}[shift={(-11.5,0)}]
\foreach \x/\y/\z in {4/-1/1,4/1/2,-4/1/3,-4/-1/4}{ \node[extra] (\z) at (\x,\y){};} \foreach \x/\y in
{1/2,2/3,3/4,4/1}{ \path[edge] (\x) -- (\y);}

\node[ver] () at (0,0.4){\mbox{{$M$} is neighborly}};

\node[ver] () at (0,-0.4){\mbox{and in {$\mathcal{K}(d)$}}};
\end{scope}

\begin{scope}[shift={(-0.5,0)}]
\foreach \x/\y/\z in {4/-1/1,4/1/2,-4/1/3,-4/-1/4}{ \node[extra] (\z) at (\x,\y){};} \foreach \x/\y in
{1/2,2/3,3/4,4/1}{ \path[edge] (\x) -- (\y);}

\node[ver] () at (0,0.4){\mbox{{$M$} is neighborly}};

\node[ver] () at (0,-0.4){\mbox{ and in {$\mathcal{H}^{d+1}$}}};
\end{scope}

\begin{scope}[shift={(-21.5,0)}]
\foreach \x/\y/\z in {3/-1/1,3/1/2,-3/1/3,-3/-1/4}{ \node[extra] (\z) at (\x,\y){};} \foreach \x/\y in
{1/2,2/3,3/4,4/1}{ \path[edge] (\x) -- (\y);}

\node[ver] () at (0,0.4){\mbox{{$M$} is }}; \node[ver] () at (0,-0.4){\mbox{tight-neighborly}};
\end{scope}

\begin{scope}[shift={(-2,-6)}]
\foreach \x/\y/\z in {5.5/-1/1,5.5/1/2,-5.5/1/3,-5.5/-1/4}{ \node[extra] (\z) at (\x,\y){};} \foreach \x/\y in
{1/2,2/3,3/4,4/1}{ \path[edge] (\x) -- (\y);}

\node[ver] () at (0,0.4){\mbox{{$M$} is {$\mathbb{F}$}-tight \& {$\beta_i(M)$} }};

\node[ver] () at (0,-0.4){\mbox{{$=0$}  for {$1 < i < d-1$}}};
\end{scope}

\begin{scope}[shift={(-21.3,-6)}]
\foreach \x/\y/\z in {3.3/-1/1,3.3/1/2,-3.3/1/3,-3.3/-1/4}{ \node[extra] (\z) at (\x,\y){};} \foreach \x/\y in
{1/2,2/3,3/4,4/1}{ \path[edge] (\x) -- (\y);}

\node[ver] () at (0,0.4){\mbox{{$M$} is strongly}};

\node[ver] () at (0,-0.4){\mbox{ minimal}};
\end{scope}

\begin{scope}[]

\path[edge] (-7.5,-2.5) -- (-5.5,-4.5);

\node[ver] () at (-7.4,-3.7){\mbox{{$\cite{dm14}$}}};

\node[ver] () at (-7.3,-2.7){\mbox{{$\nwarrow$}}};

\path[edge] (-7.5,2.5) -- (-5.5,4.5);

\node[ver] () at (-7.2,3.7){\mbox{ def.}};

\node[ver] () at (-7.3,2.7){\mbox{\boldmath{$\swarrow$}}};

\path[edge] (-7,0.5) -- (-5,0.5); \path[edge] (-7,-0.5) -- (-5,-0.5);

\path[edge] (-16,0.5) -- (-18,0.5); \path[edge] (-18,-0.5) -- (-16,-0.5);

\path[edge] (0,-2) -- (0,-4); \path[edge] (-22,-2) -- (-22,-4);

\path[edge] (-0.5,1.5) -- (-0.5,3.5); \path[edge] (0.5,1.5) -- (0.5,3.5);

\node[ver] () at (-6.85,0.47){\mbox{\boldmath{$\leftarrow$}}};

\node[ver] () at (-5.15,-0.53){\mbox{\boldmath{$\rightarrow$}}};
 \node[ver] () at (-17.75,0.47){\mbox{\boldmath{$\leftarrow$}}};
 \node[ver] () at (-16.25,-0.53){\mbox{\boldmath{$\rightarrow$}}};
 \node[ver] () at (0,-3.8){\mbox{\boldmath{$\downarrow$}}};
 \node[ver] () at (-22,-3.8){\mbox{\boldmath{$\downarrow$}}};
 \node[ver] () at (-0.5,3.3){\mbox{\boldmath{$\uparrow$}}};
 \node[ver] () at  (0.5,1.7){\mbox{\boldmath{$\downarrow$}}};
 \node[ver] () at (-6,1.1){\mbox{ def.}};
 \node[ver] () at (-6,-1.3){\mbox{ {$\cite{ka87}$}}};
 \node[ver] () at (-17.1,1.1){\mbox{ {$\cite{ns09}$}}};
 \node[ver] () at (-17.1,-1.3){\mbox{ {$\cite{ns09}$}}};
 \node[ver] () at (-1.8,2.7){\mbox{ def.}};
 \node[ver] () at (1.1,2.5){\mbox{ {$\cite{dm14}$}}};
 \node[ver] () at (1.5,-3){\mbox{ {$\cite{ef11, bd16}$}}};
 \node[ver] () at (-20.9,-3){\mbox{{$\cite{ns09}$}}};
\end{scope}
\end{tikzpicture}

\end{figure}

\item[R10.] Let $M$ be a neighbourly stacked triangulated manifold with boundary of dimension $d\geq 3$. Then $M$
is $\FF$-tight for any field $\FF$ \cite{ba14b}.

\item[R11.] For any field $\FF$ and $d\geq 4$, there does not exist a 
tight-neighborly triangulated $d$-manifold
$M$ with $\beta_1(M; \FF) =2$ \cite{si14a}.

\item[R12.] If $X$ is a triangulated 4-manifold then $f_0(X)(f_0(X)-11) \geq -15\chi(X)$. Moreover,
$f_0(X)(f_0(X)-11) = - 15\chi(X)$ implies $X \in {\cal K}(4)$ and neighborly \cite{ku95}.

\item[R13.] For $d\geq 4$, the map $M \mapsto \partial M$ is a bijection between $\overline{\mathcal K}(d+1)$ and
${\mathcal K}(d)$ \cite{bd17}.

\item[R14.] Let $G$ be a finite graph and ${\mathcal T} = \{T_i\}_{i=1}^{n}$ be a family of $(n-d-1)$-vertex
induced subtrees of $G$. Suppose that (i) any two of the $T_i$'s intersect, (ii) each vertex of $G$ is in exactly
$d+2$ members of $\mathcal T$ and (iii) for any two vertices $u\neq v$ of $G$, $u$ and $v$ are together in
exactly $d+1$ members of $\cal T$ if and only if $uv$ is an edge of $G$. Then the pure $(d+ 1)$-dimensional
simplicial complex $M$ whose facets are $\hat{u}:= \{i\in\{1, \dots, n\} \, : \, u\in T_i\}$, where $u \in V(G)$,
is an $n$-vertex neighborly member of $\overline{\mathcal K}(d+1)$, with $\Lambda(M)\cong G$ \cite{ds13}.

\item[R15.] For $k\geq 1$, let $M$ be a $(k+1)$-neighbourly and $\mathbb{F}$-orientable closed triangulated
manifold of dimension $2k$. Then $M$ is $\mathbb{F}$-tight \cite{ku95}.

\end{description}

\section{Examples}

\subsection{Some tight triangulated manifolds in {\boldmath $\mathcal{K}(d)$}}

\begin{description}
\item[E01.] For $n \geq 4$, there exist $n$-vertex neighborly orientable triangulated closed 2- \hspace{-3mm}
manifolds if and only if $n\equiv 0, 3, 4$ or 7 (mod 12) (\cite{ri74}). All are $\FF$-tight, for any field
$\FF$, by R02.

\item[E02.] For $n \geq 6$, there exist $n$-vertex neighborly non-orientable triangulated closed 2-manifolds if
and only if $n \equiv 0$ or 1 (mod 3), except for $n = 7$ (\cite{ri74}). All are $(\ZZ/2\ZZ)$-tight by R02.

\item[E03.] There is a $15$-vertex non-orientable triangulated $4$-manifold $N^{\hspace{.15mm}4}_{15}$ with
$\chi(N^{\hspace{.15mm}4}_{15}) = -4$ (\cite{bd10}). By R12, $N^{\hspace{.15mm}4}_{15}$ is neighborly and in
${\cal K}(4)$. So, $N^{\hspace{.15mm}4}_{15}$ is $(\mathbb{Z}/2\ZZ)$-tight and triangulates $(\TPSS)^{\#3}$.
Also, by R09, $N^{\hspace{.15mm}4}_{15}$ is tight-neighborly, stacked and in $\mathcal{H}^5$. \newline [There are
exactly 10 such non-orientable manifold with non-trivial $(\ZZ/3\ZZ)$- \hspace{-3mm} symmetry (\cite{si14b}).]

\item[E04.] There is a $15$-vertex orientable triangulated $4$-manifold $M^{\hspace{.15mm}4}_{15}$ with
$\chi(M^{\hspace{.15mm}4}_{15}) = -4$ (\cite{si14b}). By R12, $M^{\hspace{.15mm}4}_{15}$ is neighborly and in
${\cal K}(4)$. So, for any field $\mathbb{F}$, $M^{\hspace{.15mm}4}_{15}$ is $\mathbb{F}$-tight and triangulates
$(S^{\hspace{.15mm}3}\times S^1)^{\#3}$. Also, by R09, $M^{\hspace{.15mm}4}_{15}$ is tight-neighborly, stacked
and in $\mathcal{H}^5$. \newline [There are exactly 2 such orientable manifold with non-trivial
$(\ZZ/3\ZZ)$-symmetry (\cite{si14b}).]

\item[E05.] For $n \in \{21, 41\}$, there is an orientable $n$-vertex triangulated $4$-manifold
$M^{\hspace{.15mm}4}_{n}$ with $15\chi(M^{\hspace{.15mm}4}_{n}) = -n(n-11)$ (\cite{ds12}). By R12,
$M^{\hspace{.15mm}4}_{n}$ is neighborly and in ${\mathcal K}(4)$. By R09, $M^{\hspace{.15mm}4}_{n}$ is
$\FF$-tight, for any field $\FF$, and triangulates $(S^{\hspace{.15mm}3}\times S^1)^{\#\beta_1}$, where $\beta_1
= 1- \frac{1}{2}\chi(M^{\hspace{.15mm}4}_{n}) = \frac{1}{15}\binom{n-5}{2}$. Also, by R09,
$M^{\hspace{.15mm}4}_{n}$ is tight-neighborly, stacked and in $\mathcal{H}^5$.

\item[E06.] For $n \in \{21, 26\}$, there is a non-orientable $n$-vertex triangulated $4$-manifold
$N^{\hspace{.15mm}4}_{n}$ with $15\chi(N^{\hspace{.15mm}4}_{n}) = -n(n-11)$ (\cite{ds12}). By R12,
$N^{\hspace{.15mm}4}_{n}$ is neighborly and in ${\mathcal K}(4)$. By R09, $N^{\hspace{.15mm}4}_{n}$ is
$(\mathbb{Z}/2\ZZ)$-tight and triangulates $(\TPSS)^{\#\beta_1}$, where $\beta_1 = \frac{1}{15}\binom{n-5}{2}$.
Also, by R09, $M^{\hspace{.15mm}4}_{n}$ is tight-neighborly, stacked and in $\mathcal{H}^5$.

\item[E07.] Recently, we found 75 neighbourly, non-orientable, stacked, closed 3-manifolds. By R07, all are
tight-neighborly and $(\ZZ/2\ZZ)$-tight. These are constructed with the help of computer and R14. Among these seventy five
examples, six are 29-vertex, one is 49-vertex, fifteen are 49-vertex, forty one  are 89-vertex and twelve are 109-vertex examples \cite{bdss2}.

\end{description}

\subsection{Some infinite series of neighborly manifolds in {\boldmath $\mathcal{K}(d)$}}

\begin{description}

\item[E08.] The standard $d$-sphere $S^{\hspace{.2mm}d}_{d+2}$ (the boundary complex of a $(d+1)$-simplex) is
neighborly and in $\mathcal{H}^{d+1}$ and $\mathbb{F}$-tight for $d\geq 1$ and for any field $\mathbb{F}$.

\item[E09.] For $d\geq 2$, let $B$ be the stacked $(d+1)$-ball whose facets are $\{i, \dots, i+d+1\}$, $1\leq i
\leq 2d+3$. Let $S := \partial B$. Then $S$ is a stacked $d$-sphere. \newline Let $K^{\hspace{.2mm}d}_{2d+3}$ be obtained from $S\setminus \{\{1, \dots, d+1\}, \{2d+4, \dots, 3d+4\}\}$ by
identifying $j$ with $j+2d+3$, $1\leq j\leq d+1$. Then $K^{\hspace{.2mm}d}_{2d+3}\in {\mathcal K}(d)$ and
neighborly. Also $K^{\hspace{.2mm}d}_{2d+3}$ is orientable for $d$ even and non-orientable for $d$ odd. If $d\geq
3$ then $\beta_1(K^{\hspace{.2mm}d}_{2d+3}; \FF)=1$ for any field $\FF$ (\cite{ku86}). By R07 and R09,
$K^{\hspace{.2mm}d}_{2d+3}$ is tight-neighborly for $d\geq 3$.  If $d\geq 2$ is even then, by R02 and R09,
$K^{\hspace{.2mm}d}_{2d+3}$ is $\mathbb{F}$-tight for any field $\FF$. If $d\geq 3$ is odd then, by R07 and R09,
$K^{\hspace{.2mm}d}_{2d+3}$ is $(\ZZ/2\ZZ)$-tight. \newline [For $d\geq 3$, $K^{\hspace{.2mm}d}_{2d+3}$ is the
unique non-simply connected $(2d+3)$-vertex triangulated $d$-manifold (\cite{bd8}).]

\item[E10.] For $d\geq 2$ and $n=d^{\hspace{.3mm}2}+5d+5$, there exists an $n$-vertex  neighborly stacked
$(d+1)$-manifold with boundary $\overline{M}^{d+1}_{n}$ (in $\overline{\mathcal{K}}(5)$). Also,
$\beta_1(\overline{M}^{d+1}_{n}; \FF) = \binom{n-d-1}{2}/\binom{d +2}{2}=  d^{\hspace{.3mm}2}+5d+6$ for any field
$\FF$ (\cite{ds13}). This was constructed using R14. The dual graph $\Lambda(\overline{M}^5_{41})$ of
$\overline{M}^5_{41}$ is given in Figure 1. (For a vertex $i$ of $\overline{M}^5_{41}$, let $V_i$ be the set of
facets of $\overline{M}^5_{41}$ containing $i$. Let $T_i$ be the induced subgraph
$\Lambda(\overline{M}^5_{41})[V_i]$. Then $T_i$ is a tree.) By R10, $\overline{M}^{d+1}_{n}$ is $\FF$-tight for
any field $\FF$. The cyclic group $\ZZ/n\ZZ$ acts on $\overline{M}^{d+1}_{n}$ vertex-transitively. 


\begin{figure}[ht]
\centering
\begin{tikzpicture}[gray,thin,scale=1.4]
\draw (0,0) circle (3.0cm); \foreach \i in {0,...,40} { \draw ({(360/41*\i)}:3.0cm) node {$\bullet$}; \draw
({(360/41*\i)}:2.0cm) node {$\bullet$}; \draw ({(360/41*\i)}:2.0cm) -- ({(360/41*\i)}:3.0cm); \draw
({(360*7/41*\i)}:2.0cm) -- ({(360*7/41 *(\i+1))}:2.0cm); \draw ({(360/41*\i)}:2.2cm) node {$\circ$}; \draw
({(360/41*\i)}:2.4cm) node {$\circ$}; \draw ({(360/41*\i)}:2.6cm) node {$\circ$}; \draw ({(360/41*\i)}:2.8cm)
node {$\circ$}; }

    \foreach \i in {0} {
      \draw ({(360/41*\i)}:3.0cm) node {\color{black}{$\bullet$}}
            ({(360/41*\i)}:2.8cm) node {\color{black}{$\circ$}}
             ({(360/41*\i)}:2.6cm) node {\color{black}{$\circ$}}
             ({(360/41*\i)}:2.4cm) node {\color{black}{$\circ$}}
             ({(360/41*\i)}:2.2cm) node {\color{black}{$\circ$}}
            ({(360/41*\i)}:2.0cm) node {\color{black}{$\bullet$}}
            ({(360/41*(\i+1))}:3.0cm) node {\color{black}{$\bullet$}}
            ({(360/41*(\i+2))}:3.0cm) node {\color{black}{$\bullet$}}
            ({(360/41*(\i+2))}:2.8cm) node {\color{black}{$\circ$}}
            ({(360/41*(\i+2))}:2.6cm) node {\color{black}{$\circ$}}
            ({(360/41*(\i+2))}:2.4cm) node {\color{black}{$\circ$}}
            ({(360/41*(\i+2))}:2.2cm) node {\color{black}{$\circ$}}
            ({(360/41*(\i+3))}:3.0cm) node {\color{black}{$\bullet$}}
            ({(360/41*(\i+3))}:2.8cm) node {\color{black}{$\circ$}}
            ({(360/41*(\i+3))}:2.6cm) node {\color{black}{$\circ$}}
            ({(360/41*(\i+3))}:2.4cm) node {\color{black}{$\circ$}}
            ({(360/41*(\i+4))}:3.0cm) node {\color{black}{$\bullet$}}
            ({(360/41*(\i+4))}:2.8cm) node {\color{black}{$\circ$}}
            ({(360/41*(\i+4))}:2.6cm) node {\color{black}{$\circ$}}
            ({(360/41*(\i+5))}:3.0cm) node {\color{black}{$\bullet$}}
            ({(360/41*(\i+5))}:2.8cm) node {\color{black}{$\circ$}};

    \draw ({(360/41*(\i+7))}:2.0cm) node {\color{black}{$\bullet$}}
          ({(360/41*(\i+14))}:2.0cm) node {\color{black}{$\bullet$}}
          ({(360/41*(\i+14))}:2.2cm) node {\color{black}{$\circ$}}
          ({(360/41*(\i+21))}:2.0cm) node {\color{black}{$\bullet$}}
          ({(360/41*(\i+21))}:2.2cm) node {\color{black}{$\circ$}}
          ({(360/41*(\i+21))}:2.4cm) node {\color{black}{$\circ$}}
          ({(360/41*(\i+28))}:2.0cm) node {\color{black}{$\bullet$}}
          ({(360/41*(\i+28))}:2.2cm) node {\color{black}{$\circ$}}
          ({(360/41*(\i+28))}:2.4cm) node {\color{black}{$\circ$}}
          ({(360/41*(\i+28))}:2.6cm) node {\color{black}{$\circ$}}
          ({(360/41*(\i+35))}:2.0cm) node {\color{black}{$\bullet$}}
          ({(360/41*(\i+35))}:2.2cm) node {\color{black}{$\circ$}}
          ({(360/41*(\i+35))}:2.4cm) node {\color{black}{$\circ$}}
          ({(360/41*(\i+35))}:2.6cm) node {\color{black}{$\circ$}}
          ({(360/41*(\i+35))}:2.8cm) node {\color{black}{$\circ$}};

    }

    \foreach \i in {0} {
        \draw [black,thick] ({(360/41*\i)}:3.0cm) arc ({360/41*\i}:{360/41*(\i+5)}:3.0cm)
        ({(360/41*\i)}:3.0cm) -- ({(360/41*\i)}:2.0cm)
        -- ({360/41*(\i+7)}:2.0cm)
        -- ({360/41*(\i+14)}:2.0cm)
        -- ({360/41*(\i+21)}:2.0cm)
        -- ({360/41*(\i+28)}:2.0cm)
        -- ({360/41*(\i+35)}:2.0cm)
        ({360/41*(\i+2)}:3.0cm) -- ({360/41*(\i+2)}:2.2cm)
        ({360/41*(\i+3)}:3.0cm) -- ({360/41*(\i+3)}:2.4cm)
        ({360/41*(\i+4)}:3.0cm) -- ({360/41*(\i+4)}:2.6cm)
        ({360/41*(\i+5)}:3.0cm) -- ({360/41*(\i+5)}:2.8cm)
        ({360/41*(\i+14)}:2.0cm) -- ({360/41*(\i+14)}:2.2cm)
        ({360/41*(\i+21)}:2.0cm) -- ({360/41*(\i+21)}:2.4cm)
        ({360/41*(\i+28)}:2.0cm) -- ({360/41*(\i+28)}:2.6cm)
        ({360/41*(\i+35)}:2.0cm) -- ({360/41*(\i+35)}:2.8cm);
    }
\end{tikzpicture}
\caption{Graph $\Lambda(\overline{M}^{5}_{41})$ and the tree $T_1$ in \textcolor{black}{black}} 
\end{figure}
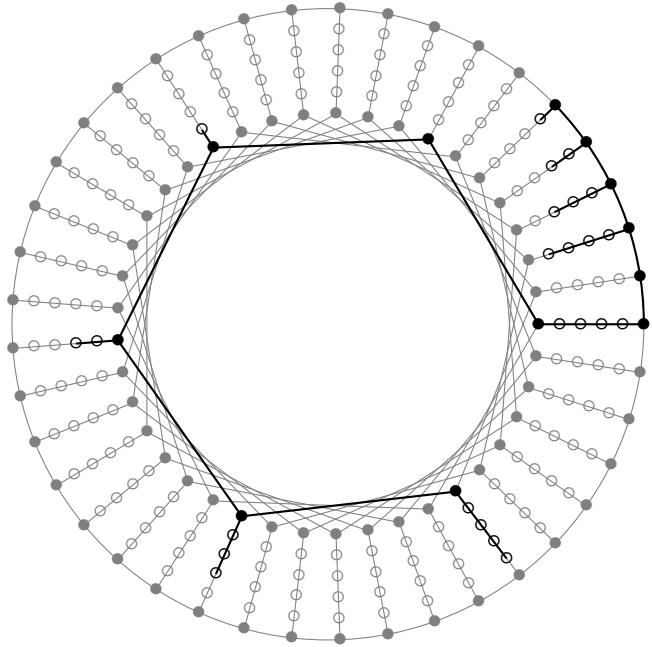


\item[E11.] For $d\geq 2$ and $n=d^{\hspace{.3mm}2}+5d+5$, let $M^{d}_{n} :=\partial \overline{M}^{d+1}_{n}$,
where $\overline{M}^{d+1}_{n}$ as in E10. Then $M^{d}_{n} \in \mathcal{H}^{d+1}$ and neighborly. Moreover,
$M^{d}_{n}$ is orientable for $d$ even and non-orientable for $d$ odd (\cite{ds13}). By R02, R07 and R09, we have
the following: (a) $M^d_{n}$ is $\FF$-tight for any field $\FF$ if $d$ is even and is $(\ZZ/2\ZZ)$-tight if $d$
is odd. (b) $M^d_{n}$ is tight-neighborly for $d\geq 3$. (c) $\beta_1(M^d_{n}; \FF) =\binom{n-d-1}{2}/\binom{d
+2}{2}= d^{\hspace{.3mm}2}+5d+6$ for $d\geq 3$ and any field $\FF$. (d) If $d\geq 2$ is even then $M^d_{n}$
triangulates $(S^{\hspace{.1mm}d -1}\!\times S^1)^{\# \beta}$ and if $d\geq 3$ is odd then $M^d_{n}$ triangulates
$(\TPSSD)^{\#\beta}$, where $\beta= d^{\hspace{.3mm}2}+5d+6$. (e) $M^d_{n}$ is strongly minimal for all $d$. (f)
$\mathbb{Z}/n\ZZ$ acts vertex-transitively on $M^d_{n}$ for all $d$.

\end{description}

\newpage

\subsection{Known neighborly members of {\boldmath${\mathcal K}(d)$} 
for {\boldmath $d\geq 3$}}

\begin{table}[htbp]
\begin{center}

\begin{tabular}{|c|c|c|c|c|l|}
\hline
&&&&&\\[-3mm]
$\beta_1(K)$ & $d$ & $f_0(K)$ & $K$ & $|K|$ & References\\[0.5mm]
\hline
&&&&&\\[-3mm]
0 & $d$ & $d+2$ & {$S^{\hspace{.1mm}d}_{d+2}$} &
\boldmath{$S^{\hspace{.1mm}d}$} & E08 \\[0.5mm]
\hline
&&&&&\\[-3mm]
1 & $d$ even & $2d+3$ & {$K^{d}_{2d+3}$} &
$S^{d-1}\times S^1$ & E09, \cite{ku86} \\[1mm]
,, & $d$ odd & ,, & ,, & {$S^{d-1}\simtimes S^1$} & ,, \\[0.5mm]
\hline
&&&&&\\[-2.5mm]
$\frac{1}{10}\binom{f_0-4}{2}$ & 3 & $20\ell +9$, $1\leq \ell \leq 5$ & 
$\exists$ 75 & $(S^2\simtimes S^1)^{\#\beta_1}$ &  E07, \cite{bdss2} 
\\[1mm]
\hline
&&&&&\\[-3mm]
2 & $d\geq 4$ & $-$  & $\not\!\exists$ & & R11, \cite{si14a} \\[0.5mm]
\hline
&&&&&\\[-3mm]
3 & 4 & 15 & {$N^4_{15}$} & {$(S^3 \simtimes S^1)^{\#3}$} & E03, \cite{bd10} \\[1mm]
,, & ,, & ,, & {$M^4_{15}$} & {$(S^3\times S^1)^{\#3}$} & E04, \cite{si14b} \\[0.5mm]
\hline
&&&&&\\[-3mm]
8 & 4 & 21 & {$M^4_{21}$} & {$(S^3\times S^1)^{\#8}$} &E05, \cite{ds12} \\[1mm]
,, & ,, & ,, & {$N^4_{21}$} & {$(S^3 \simtimes S^1)^{\#8}$} & E06, ,, \\[0.5mm]
\hline
&&&&&\\[-3mm]
14 & 4 & 26 & {$N^4_{26}$} & {$(S^3 \simtimes S^1)^{\#14}$} & E06, \cite{ds12} \\[0.5mm]
\hline
&&&&&\\[-3mm]
 $2\binom{d+3}{2}$ & $d$ even & $d^2+5d+5$ & {$M^{d}_{n}$} &
{$(S^{d-1}\times S^1)^{\#\beta_1}$}  & E11, \cite{ds13}\\[1mm]
 ,, & $d$ odd & ,, & ,, & {$(S^{d-1}\simtimes S^1)^{\#\beta_1}$} & ,,\\[0.5mm]
\hline
\end{tabular}

\label{tbl:tbl1}
\end{center}
\end{table}

\subsection{Some tight triangulated manifolds outside {\boldmath $\mathcal{K}(d)$}}

\begin{description}

\item[E12.] Lutz constructed two 12-vertex triangulations of $S^{\,2} \times S^{\,3}$ \cite{lu1}. Both are
$\FF$-tight for any field $\FF$.


\item[E13.] Only finitely many $2k$-dimensional $(k+ 1)$-neighbourly triangulated closed manifolds are known for
$k \geq 2$. By R15, they are all $\FF$-tight. These examples are\,:


 \begin{enumerate}[{(a)}]

 \item the 9-vertex triangulation of $\mathbb{C}\mathbb{P}^{\,2}$ due to K\"{u}hnel \cite{kb},


 \item six 15-vertex triangulations of homology $\mathbb{H P}^{\,2}$ (three due to Brehm and
K\"{u}hnel \cite{bk2} and three more due to Lutz \cite{lu2}),


 \item the 16-vertex triangulation of a $K3$-surface due to Casella and K\"{u}hnel \cite{ck}, and


 \item two 13-vertex triangulations of $S^{\,3} \times S^{\,3}$ due to Lutz \cite{lu1}.
\end{enumerate}


\item[E14.] Apart from above, we know only two tight triangulated manifolds. These are\,:


\begin{enumerate}[{(a)}]

\item A 15-vertex triangulation of $(\TPSS) \# (\mathbb{C}\mathbb{P}^{\,2})^{\#5}$ due to Lutz \cite{lu1}. It is
non-orientable and $(\mathbb{Z}/2\ZZ)$-tight.


\item A 13-vertex triangulation of $SU(3)/SO(3)$ due to Lutz \cite{lu1}. It is 3-neighbourly, orientable and
$(\mathbb{Z}/2\ZZ)$-tight.
\end{enumerate}

\end{description}

\bigskip

\noindent {\bf Acknowledgements\,:} The author is supported by DIICCSRTE, Australia and DST, India, under the
Australia-India Strategic Research Fund (project AISRF06660) and by the UGC Centre for Advanced Studies.

{\small

}

\end{document}